\documentclass[12pt]{article}
\usepackage{amssymb}

\newcommand{\be}{\begin{equation}}
\newcommand{\ee}{\end{equation}}
\newcommand{\bea}{\begin{eqnarray}}
\newcommand{\eea}{\end{eqnarray}}
\newcommand{\barray}{\begin{array}}
\newcommand{\earray}{\end{array}}
\newcommand{\pa}{\partial}
\newcommand{\nn}{\nonumber}
\newcommand{\bitem}{\begin{itemize}}
\newcommand{\eitem}{\end{itemize}}
\newtheorem{teo}{Theorem}[section]
\newcommand{\bt}{\begin{teo}}
\newcommand{\et}{\end{teo}}
\newtheorem{Def}{Definition}[section]
\newcommand{\bd}{\begin{Def}}
\newcommand{\ed}{\end{Def}}
\newtheorem{lem}{Lemma}[section]
\newcommand{\bl}{\begin{lem}}
\newcommand{\el}{\end{lem}}
\newtheorem{prop}{Proposition}[section]
\newcommand{\bp}{\begin{prop}}
\newcommand{\ep}{\end{prop}}
\newtheorem{cor}{Corollary}[section]
\newcommand{\bc}{\begin{cor}}
\newcommand{\ec}{\end{cor}}
\newtheorem{ex}{Example}[section]
\newcommand{\bex}{\begin{ex}}
\newcommand{\eex}{\end{ex}}
\newtheorem{rem}{Remark}[section]
\newcommand{\br}{\begin{rem}}
\newcommand{\er}{\end{rem}}

\begin{document}

\begin{center}
{\Large \textbf{Realization of Frobenius manifolds \\ as
submanifolds in pseudo-Euclidean spaces\footnotetext[1]{The work was
supported by the Max-Planck-Institut f\"ur Mathematik (Bonn,
Germany), by the Russian Foundation for Basic Research (project
no.~08-01-00054) and by a grant of the President of the Russian
Federation (project no. NSh-1824.2008.1).}}}
\end{center}

\smallskip

\begin{center}
{\large {O. I. Mokhov}}
\end{center}

\smallskip

\begin{abstract}
We introduce a class of $k$-potential submanifolds in
pseudo-Euclidean spaces and prove that for an arbitrary positive
integer $k$ and an arbitrary nonnegative integer $p$, each
$N$-dimensional Frobenius manifold can always be locally realized as
an $N$-dimensional $k$-potential submanifold in $((k + 1) N +
p)$-dimensional pseudo-Euclidean spaces of certain signatures. For
$k = 1$ this construction was proposed by the present author in a
previous paper (2006). The realization of concrete Frobenius
manifolds is reduced to solving a consistent linear system of
second-order partial differential equations.
\end{abstract}

\smallskip

\begin{flushright}
to Vladimir Igorevich Arnold
\end{flushright}

\smallskip

\smallskip

\section{Introduction}

In this paper we develop and significantly generalize the
construction proposed earlier by the present author in \cite{1} and
generated by a deep nontrivial relationship discovered in \cite{1}
between the theory of Frobenius manifolds, the associativity
equations of two-dimensional topological quantum field theories (the
Witten--Dijkgraaf--Verlinde--Ver\-linde equations), and the
Dubrovin--Frobenius structures on the one hand and the theory of
submanifolds in pseudo-Euclidean spaces on the other hand. In this
connection we construct new very natural integrable {\it
$k$-potential\/} reductions of the fundamental Gauss--Codazzi--Ricci
equations and new interesting integrable classes of submanifolds in
pseudo-Euclidean spaces. These classes are important for
applications. In particular, in this paper we introduce a new
integrable class of {\it $k$-potential\/} submanifolds in
pseudo-Euclidean spaces and prove that for an arbitrary positive
integer $k$ and an arbitrary nonnegative integer $p$, each
$N$-dimensional Frobenius manifold can always be locally realized as
an $N$-dimensional {\it $k$-potential\/} submanifold in $((k + 1) N
+ p)$-dimensional pseudo-Euclidean spaces of certain signatures. For
$k = 1$ this construction was proposed by the present author in
\cite{1}. The realization of any concrete $N$-dimensional Frobenius
manifold as an $N$-dimensional $k$-potential submanifold in $((k +
1) N + p)$-dimensional pseudo-Euclidean spaces is reduced to solving
a consistent linear system of second-order partial differential
equations, which can be solved explicitly in elementary and special
functions. First of all, we prove that for an arbitrary positive
integer $k$ and an arbitrary nonnegative integer $p$, the
associativity equations of two-dimensional topological quantum field
theories (the Witten--Dijkgraaf--Verlinde--Verlinde (WDVV)
equations, see \cite{1d}--\cite{4d}) for a function (a {\it
potential}\,) $\Phi = \Phi (u^1, \ldots, u^N)$, \be \sum_{k = 1}^N
\sum_{l = 1}^N {\pa^3 \Phi \over \pa u^i \pa u^j \pa u^k} \eta^{kl}
{\pa^3 \Phi \over \pa u^l \pa u^m \pa u^n} = \sum_{k = 1}^N \sum_{l
= 1}^N {\pa^3 \Phi \over \pa u^i \pa u^m \pa u^k} \eta^{kl} {\pa^3
\Phi \over \pa u^l \pa u^j \pa u^n}, \label{ass1} \ee where
$\eta^{ij}$ is an arbitrary constant nondegenerate symmetric matrix,
$\eta^{ij} = {\rm const}$, $\det \, (\eta^{ij}) \neq 0$, $\eta^{ij}
= \eta^{ji}$, are very natural integrable $k$-potential reductions
of the fundamental nonlinear equations in the submanifold theory
(the corresponding Gauss--Codazzi--Ricci equations) that describe
$N$-dimensional submanifolds in $((k + 1) N + p)$-dimensional
pseudo-Euclidean spaces. The WDVV equations give a natural and
important nontrivial integrable class of {\it $k$-potential\/}
$N$-dimensional submanifolds of codimension $k N + p$ in $((k + 1) N
+ p)$-dimensional pseudo-Euclidean spaces. For the special case $k =
1$ and $p = 0$ all these statements were formulated and proved in
\cite{1} and \cite{dga}, where such $N$-dimensional submanifolds of
codimension $N$ were called {\it potential\/}. All $k$-potential
submanifolds have natural differential-geometric special structures
of Frobenius algeb\-ras (the Dubrovin--Frobenius structures) on
their tangent spaces. These Dubrovin--Frobenius structures are
generated by the corresponding flat first fundamental forms and some
sets of the second fundamental forms of the submanifolds (the
structural constants of a Frobenius algebra are given and duplicated
by sets of the Weingarten operators of the $k$-potential
submanifold).

We recall that each solution $\Phi (u^1, \ldots, u^N)$ of the
associativity equations (\ref{ass1}) gives Dubrovin--Frobenius
structures, i.e., specific $N$-parameter deformations of Frobenius
algebras; in our case these algebras are commutative associative
algebras equipped with nondegenerate invariant symmetric bilinear
forms. Indeed, consider the algebras $A (u)$ in an $N$-dimensional
vector space with the basis $e_1, \ldots, e_N$ and multiplication
(see \cite{1d}) \be e_i \circ e_j = c^k_{ij} (u) e_k, \ \ \ \
c^k_{ij} (u) = \eta^{ks} {\pa^3 \Phi \over \pa u^s \pa u^i \pa u^j}.
\label{01} \ee For all values of the parameters $u = (u^1, \ldots,
u^N)$, the algebras $A (u)$ are commutative, $e_i \circ e_j = e_j
\circ e_i,$ and the associativity condition \be (e_i \circ e_j)
\circ e_k = e_i \circ (e_j \circ e_k) \label{02} \ee in the algebras
$A (u)$ is equivalent to the WDVV equations (\ref{ass1}). The
inverse $\eta_{ij}$ of the matrix $\eta^{ij}$, $\eta^{is} \eta_{sj}
= \delta^i_j$, defines a nondegenerate invariant symmetric bilinear
form on the algebras $A (u)$, \be \langle e_i, e_j \rangle =
\eta_{ij}, \ \ \ \ \langle e_i \circ e_j, e_k \rangle = \langle e_i,
e_j \circ e_k \rangle. \label{03} \ee Recall that locally any
Frobenius manifold has a Dubrovin--Frobenius structure (see
\cite{1d}); namely, the tangent space at every point $u = (u^1,
\ldots, u^N)$ of any Frobenius manifold possesses the structure of a
Frobenius algebra (\ref{01})--(\ref{03}), which is determined by a
certain solution of the associativity equations (\ref{ass1}) and
smoothly depends on the point. In this paper we prove that for an
arbitrary positive integer $k$ and an arbitrary nonnegative integer
$p$, each $N$-dimensional Frobenius manifold can always be locally
represented as an $N$-dimensional $k$-potential submanifold in some
$((k + 1) N + p)$-dimensional pseudo-Euclidean spaces. The
corresponding representations of any given Frobenius manifold are
parametrized by the set of admissible Gram matrices of the scalar
products of basic vectors in the normal spaces of the corresponding
$k$-potential submanifolds; in particular, for $p = 0$ they are
parametrized by the set of arbitrary nondegenerate symmetric
constant $k \times k$ matrices determining all the admissible $k N
\times k N$ Gram matrices of the scalar products of basic vectors in
the normal spaces of the corresponding $k$-potential submanifolds.
If for an arbitrary Frobenius manifold we fix a certain admissible
$(k N + p) \times (k N + p)$ Gram matrix of the scalar products of
basic vectors in the normal spaces of the corresponding
$k$-potential submanifolds, then the corresponding $k$-potential
submanifold realizing this Frobenius manifold and having the given
Gram matrix of basic vectors in the normal spaces is determined by
our construction uniquely up to motions in the corresponding ambient
pseudo-Euclidean space. We note that an alternative approach to the
description of the submanifolds realizing arbitrary Frobenius
manifolds is developed by the present author in \cite{m}.

\section{Frobenius algebras and
Frobenius manifolds}

\subsection{Frobenius algebras}

In the mathematical literature there are various widely spread
approaches to the notion of Frobenius algebra and different
definitions of Frobenius algebras not always requiring even
associativity of the algebra, to say nothing of the requirement of
presence of a unit in the algebra, symmetry of the corresponding
bilinear form, and commutativity of the algebra. Therefore, we give
here necessary definitions that will be used in this article. The
presence of a special nondegenerate bilinear form that is compatible
with the multiplication in the algebra is a common feature of all
definitions of Frobenius algebras.

Let us consider a finite-dimensional algebra ${\mathcal{A}}$ (with
multiplication $\circ$) over a field ${\mathbb{K}}$ (in this paper
we consider algebras only over ${\mathbb{R}}$ or ${\mathbb{C}}$).

\bd {\rm A bilinear form $f : {\mathcal{A}} \times {\mathcal{A}}
\rightarrow {\mathbb{K}}$ in an algebra ${\mathcal{A}}$ is called
{\it invariant\/} (or {\it associative\/}) if \be f (a \circ b, c) =
f (a, b \circ c) \ee for all $a, b, c \in {\mathcal{A}}$.}\ed

\bd {\rm A finite-dimensional associative algebra ${\mathcal{A}}$
over a field ${\mathbb{K}}$ that is equipped with a nondegenerate
invariant symmetric bilinear form is called {\it Frobenius\/}.}\ed

Note that we do not require the presence of a unit in a Frobenius
algebra.

\bex {\bf Matrix algebra $M_n ({\mathbb{K}})$}.

{\rm Consider the algebra $M_n ({\mathbb{K}})$ (the algebra of $n
\times n$ matrices over a field ${\mathbb{K}}\,$), the linear
functional (trace of matrices)
$$\theta (a) = {\rm Tr\,} (a), \ a \in M_n ({\mathbb{K}}),$$
and the bilinear form $f (a, b) = \theta (a b)$. This bilinear form
is invariant, since the matrix algebra is associative. It is easy to
prove that this bilinear form is nondegenerate, and the matrix
algebra $(M_n ({\mathbb{K}}), f)$ is a noncommutative Frobenius
algebra with a unit over a field ${\mathbb{K}}\,$. Note that the
bilinear form $f (a, b) = \theta (a b)$ is symmetric, $\theta (a b)
= \theta (b a)$.} \eex

\bex  {\bf Group algebra ${\mathbb{K}}G$}.

{\rm Let $G$ be a finite group. Consider the group algebra
${\mathbb{K}}G$ over a field ${\mathbb{K}}\,$, $${\mathbb{K}}G = \{a
\ | \ a = \sum_{g \in G} \alpha_g g, \  \alpha_g \in {\mathbb{K}}
\}.$$ Obviously, ${\mathbb{K}}G$ is an associative algebra with unit
over the field ${\mathbb{K}}\,$. Let $e$ be the unit of the group
$G$. Consider the linear functional $$ \theta (a) = \alpha_e (a), \
\ a = \sum_{g \in G} \alpha_g (a) g \in {\mathbb{K}}G, \ \ \alpha_g
(a) \in {\mathbb{K}},$$ and the bilinear form $f (a, b) = \theta (a
b)$. This bilinear form is invariant, since the group algebra is
associative. It is easy to prove that this bilinear form is
nondegenerate. Indeed, we have $$f (g^{-1}, a) = \theta (g^{-1} a) =
\alpha_g (a)$$ for all $g \in G$. Therefore, if $f (g, a) = \theta
(g a) = 0$ for all $g \in G$, then $\alpha_g (a) = 0$ for all $g \in
G$, i.e., $a = 0$. Hence the bilinear form $f$ is nondegenerate, and
the group algebra $({\mathbb{K}}G, f)$ is a noncommutative Frobenius
algebra with a unit over a field ${\mathbb{K}}\,$ (it is commutative
only for Abelian groups $G$). Note that the bilinear form $f (a, b)
= \theta (a b)$ is symmetric for any group $G$, $\theta (a b) =
\theta (b a)$.} \eex

\subsection{Frobenius manifolds}

It would be quite natural to call a manifold {\it Frobenius} if each
tangent space at any point of this manifold is equipped with a
Frobenius algebra structure that depends smoothly on the point of
the manifold. However, a remarkable and very fruitful theory of
Frobenius manifolds with very special Frobenius structures was
constructed by Dubrovin (see \cite{1d}) in connection with
two-dimensional topological quantum field theories and quantum
cohomology, and it is these manifolds that were called Frobenius.
Such Frobenius manifolds play an important role in singularity
theory, enumerative geometry, theory of Gromov--Witten invariants,
quantum cohomology theory, to\-pological quantum field theories, and
in various other domains of modern differential geometry and
mathematical and theoretical physics. In this paper we follow the
definition of \cite{1d}, but we do not impose some very severe
Dubrovin's constraints on Frobenius manifolds (in particular, we do
not require quasihomogeneity, the presence of a special Euler vector
field, and the presence of a covariantly constant unit in the
Frobenius algebra on the tangent spaces of the manifold). We will
call the corresponding structures on manifolds {\it
Dubrovin--Frobenius structures}.

\bd {\rm \cite{1d} An $N$-dimensional pseudo-Riemannian manifold
$M^N$ with a metric $g$ and an algebra structure $(T_u M, \circ)$,
$T_u M \times T_u M \stackrel{\circ}{\rightarrow} T_u M$, that is
defined on each tangent space $T_u M$ and depends smoothly on the
point $u \in M^N$ is called {\it Frobenius\/} if

(I) the pseudo-Riemannian metric $g$ is a nondegenerate invariant
symmetric bilinear form on each tangent space $T_u M$, \be g (X
\circ Y, Z) = g (X, Y \circ Z) \label{f1a} \ee for all vector fields
$X, Y,$ and $Z$ on $M^N$;

(II) the algebra $(T_u M, \circ)$ is commutative at each point $u
\in M^N$, \be X \circ Y = Y \circ X \ee for all vector fields $X$
and $Y$ on $M^N$;

(III) the algebra $(T_u M, \circ)$ is associative at each point $u
\in M^N$, \be (X \circ Y) \circ Z = X \circ (Y \circ Z) \ee for all
vector fields $X, Y,$ and $Z$ on $M^N$;

(IV) the metric $g$ is flat;

(V) $A (X, Y, Z) = g (X \circ Y, Z)$ is a symmetric tensor on $M^N$
such that the tensor $(\nabla_W A) (X, Y, Z)$ is symmetric with
respect to all vector fields $X, Y, Z,$ and $W$ on $M^N$ ($\nabla$
is the covariant differentiation generated by the Levi-Civita
connection of the metric $g$).} \ed

It is obvious that conditions (I)--(III) mean that at each point $u
\in M^N$ the algebra $(T_u M, \circ, g)$ is a commutative Frobenius
algebra.

\subsection{Associativity equations}

Let us consider an arbitrary Frobenius manifold, i.e., an arbitrary
manifold satisfying conditions (I)--(V). Let $u = (u^1, \ldots,
u^N)$ be arbitrary flat coordinates of the flat metric $g$. In any
flat local coordinates, the metric $g (u)$ is a constant
nondegenerate symmetric matrix $\eta_{ij}$, $\eta_{ij} = \eta_{ji},$
$\det (\eta_{ij}) \neq 0,$ $\eta_{ij} = \ {\rm const}$, $g (X, Y) =
\eta_{ij} X^i (u) Y^j (u)$.

In these flat local coordinates, for structural functions $c^i_{jk}
(u)$ of the Frobenius algebra on the manifold,
$$X \circ Y = W, \ \ \ W^i (u) = c^i_{jk} (u) X^j (u) Y^k (u),$$
and for the symmetric tensor $A_{ijk} (u)$, we have \bea && A (X, Y,
Z) = A_{ijk} (u) X^i (u) Y^j (u) Z^k (u)
= g (X \circ Y, Z) = \nn\\
&& = g (W, Z) = \eta_{ij} W^i (u) Z^j (u) = \eta_{ij} c^i_{kl} (u)
X^k (u) Y^l (u) Z^j (u);\nn \eea i.e., \be A_{ijk} (u) = \eta_{sk}
c^s_{ij} (u). \label{str} \ee

According to condition (V), $(\nabla_l A_{ijk}) (u)$ is a symmetric
tensor; i.e., in the flat local coordinates we also have
$${\pa A_{ijk} \over \pa u^l} = {\pa A_{ijl} \over \pa u^k}.$$

Hence there locally exists a function (a {\it potential\/}) $\Phi
(u)$ such that
$$A_{ijk} (u) = {\pa^3 \Phi \over \pa u^i \pa u^j \pa u^k}.$$

From relation (\ref{str}) for the structural functions $c^i_{jk}
(u)$ we obtain \be c^i_{jk} (u) = \eta^{is} A_{sjk} (u) = \eta^{is}
{\pa^3 \Phi \over \pa u^s \pa u^j \pa u^k}, \label{str2} \ee where
the matrix $\eta^{ij}$ is the inverse of the matrix $\eta_{ij}$,
$\eta^{is} \eta_{sj} = \delta^i_j$.

For any values of the parameters $u = (u^1, \ldots, u^N)$, the
structural functions (\ref{str2}) give a commutative algebra \be
\pa_i \circ \pa_j = c^k_{ij} (u) \pa_k = \eta^{ks} {\pa^3 \Phi \over
\pa u^s \pa u^i \pa u^j} \pa_k \label{frob} \ee equipped with a
symmetric invariant nondegenerate bilinear form \be \langle \pa_i,
\pa_j \rangle = \eta_{ij} \label{form} \ee for any constant
nondegenerate symmetric matrix $\eta_{ij}$ and for any function
$\Phi (u)$, but, generally speaking, this algebra is not
associative. All conditions (I)--(V) except the associativity
condition (III) are obviously satisfied for all these $N$-parameter
deformations of nonassociative algebras.

The associativity condition (III) is equivalent to a nontrivial
overdetermined system of nonlinear partial differential equations
for the potential $\Phi (u)$, \be \sum_{k = 1}^N \sum_{l = 1}^N
{\pa^3 \Phi \over \pa u^i \pa u^j \pa u^k} \eta^{kl} {\pa^3 \Phi
\over \pa u^l \pa u^m \pa u^n} = \sum_{k = 1}^N \sum_{l = 1}^N
{\pa^3 \Phi \over \pa u^i \pa u^m \pa u^k} \eta^{kl} {\pa^3 \Phi
\over \pa u^l \pa u^j \pa u^n},  \label{ass1a} \ee which is well
known as the system of associativity equations of two-dimensional
topological quantum field theories (the WDVV equations, see
\cite{1d}--\cite{4d}, \cite{6m}--\cite{5}). The system of
associativity equations (\ref{ass1a}) is consistent, integrable by
the inverse scattering method, and possesses a rich set of
nontrivial exact solutions (see \cite{1d}).

It is obvious that each solution $\Phi (u^1, \ldots, u^N)$ of the
associativity equations (\ref{ass1a}) gives $N$-parameter
deformations of commutative Frobenius algebras (\ref{frob}) equipped
with nondegenerate invariant symmetric bilinear forms (\ref{form}).
These Dubrovin--Frobe\-nius structures satisfy all conditions
(I)--(V).

Further in this paper we show that the associativity equations
(\ref{ass1a}) are natural reductions of the fundamental nonlinear
equations in the theory of submanifolds in pseudo-Euclidean spaces
and give a natural class of {\it $k$-potential\/} submanifolds. All
$k$-potential submanifolds in pseudo-Euclidean spaces have natural
differential-geometric structures of Frobenius algebras
(\ref{frob}), (\ref{form}) on their tangent spaces. These
differential-geometric Dubrovin--Frobenius structures are generated
by the flat first fundamental forms and the sets of the second
fundamental forms on the submanifold (the structural constants of
the Frobenius algebra are given by the Weingarten ope\-rators of the
submanifold).

A great number of concrete examples of Frobenius manifolds and
solutions of the associativity equations are given in \cite{1d}.
Consider here one simple but important example from \cite{1d}. Let
$N = 3$, let the metric $\eta_{ij}$ be antidiagonal, \be (\eta_{ij})
=
\left ( \begin{array} {ccc} 0&0&1\\
0&1&0\\
1&0&0
\end{array} \right ),
\ee and let $e_1$ be a unit in the Frobenius algebra (\ref{frob}),
(\ref{form}). In this case, the function (the potential) $\Phi (u)$
has the form
$$
\Phi (u) = {1 \over 2} (u^1)^2 u^3 + {1 \over 2} u^1 (u^2)^2 + f
(u^2, u^3),
$$ and the associativity equations
(\ref{ass1a}) for the function $\Phi (u)$ are equivalent to the
following remarkable equation for the function $f (u^2, u^3)$ (see
\cite{1d}): \be {\pa^3 f \over \pa (u^3)^3} = \left ( {\pa^3 f \over
\pa (u^2)^2 \pa u^3} \right )^2 - {\pa^3 f \over \pa (u^2)^3} {\pa^3
f \over \pa u^2 \pa (u^3)^2}. \label{f} \ee This equation is
connected with the quantum cohomology of projective plane and
classical problems of enumerative geometry (see \cite{6}). In
particular, all nontrivial polynomial quasihomogeneous solutions of
equation (\ref{f}) are described in \cite{1d}: \be f = {1 \over 4}
(u^2)^2 (u^3)^2 + {1 \over 60} (u^3)^5,\ \ \ f = {1 \over 6} (u^2)^3
u^3 + {1 \over 6} (u^2)^2 (u^3)^3 + {1 \over 210} (u^3)^7,
\label{sol1} \ee \be f = {1 \over 6} (u^2)^3 (u^3)^2 + {1 \over 20}
(u^2)^2 (u^3)^5 + {1 \over 3960} (u^3)^{11}. \label{sol2} \ee

As shown by the author in \cite{7} (see also \cite{7a}--\cite{8}),
equation (\ref{f}) is equivalent to the integrable nondiagonalizable
homogeneous system of hydrodynamic type \be \left (
\begin{array} {c} a\\ b\\ c
\end{array} \right )_{u^3} =
\left ( \begin{array} {ccc} 0 & 1 & 0\\  0 & 0 & 1\\
- c & 2b & - a
\end{array} \right )  \left ( \begin{array} {c} a\\ b\\ c
\end{array} \right )_{u^2}, \label{shdt}
\ee \be a = {\pa^3 f \over \pa (u^2)^3}, \ \ \ b = {\pa^3 f \over
\pa (u^2)^2 \pa u^3},\ \ \ c = {\pa^3 f \over \pa u^2 \pa (u^3)^2}.
\ee In this case, the Weingarten operators of potential submanifolds
that realize the corres\-ponding Frobenius manifolds have the form
\be (w_1)^i_j (u) = \delta^i_j, \ \ \ (w_2)^i_j (u) =
\left ( \begin{array} {ccc} 0 & b & c\\  1 & a & b\\
0 & 1 & 0
\end{array} \right ),\ \ \
(w_3)^i_j (u) = \left ( \begin{array} {ccc}
0 & c & b^2 - a c\\  0 & b & c\\
1 & 0 & 0
\end{array} \right ).
\ee

For concrete solutions of the associativity equation (\ref{f}), in
particular, for (\ref{sol1}) and (\ref{sol2}), or for any concrete
solutions of the system of hydrodynamic type (\ref{shdt}), the
corresponding linear systems that provide explicit realizations of
the corresponding Frobenius manifolds as $k$-potential submanifolds
in pseudo-Euclidean spaces can be solved explicitly in elementary
and special functions.

\section{General fundamental equations of the local \\ theory of totally
nonisotropic submanifolds \\ in pseudo-Eucli\-dean spaces}

Let us consider an arbitrary totally nonisotropic smooth
$N$-dimensional submanifold $M^N$ in an $(N + L)$-dimensional
pseudo-Euclidean space $\mathbb{E}^{N + L}_q$ of arbitrary signature
$q$, $M^N \subset \mathbb{E}^{N + L}_q$. Recall that a submanifold
of a pseudo-Euclidean space is called {\it totally nonisotropic} if
it is not tangent to isotropic cones of the ambient pseudo-Euclidean
space at any of its points. A submanifold of a pseudo-Euclidean
space is totally nonisotropic if and only if the metric induced on
the submanifold from the ambient pseudo-Euclidean space (the first
fundamental form of the submanifold) is nondegenerate. Note that in
this paper we consider only the local theory of submanifolds.

Let $(z^1, \ldots, z^{N + L})$ be pseudo-Euclidean coordinates in
$\mathbb{E}^{N + L}_q$, and let the subma\-nifold $M^N$ be given
locally by a smooth vector function $r (u^1, \ldots, u^N)$ of $N$
independent variables $u^1, \ldots, u^N$ (local coordinates on the
submanifold $M^N$), $r (u^1, \ldots, u^N) = (z^1 (u^1, \ldots, u^N),
\ldots, z^{N + L} (u^1, \ldots, u^N))$, ${\rm rank\,} (\partial z^i
/ \partial u^j) = N$, $1 \leq i \leq N + L$, $1 \leq j \leq N$. In
this case $\partial r /
\partial u^i = r_i (u)$, $1 \leq i \leq N$, is a basis of the tangent space
$\mathbb{T}_u$ at each point $u = (u^1, \ldots, u^N)$ of the
submanifold $M^N$, and $g_{ij} (u) = (r_i, r_j)$, $1 \leq i, j \leq
N$, is the first fundamental form of the submanifold $M^N$, where
$(\cdot, \cdot)$ is the pseudo-Euclidean scalar product in
$\mathbb{E}^{N + L}_q$. In the normal space $\mathbb{N}_u$ of the
submanifold $M^N$ at each point $u$, we fix an arbitrary basis $n_1
(u), \ldots, n_L (u)$ that depends smoothly on the point $u$.
Consider the corresponding matrix of scalar products of the basis
vectors in the normal spaces on the submanifold $M^N$ (we will also
call it the Gram matrix in the normal spaces of the submanifold
$M^N$, although in this case the scalar product is, generally
speaking, pseudo-Euclidean), i.e., the matrix of functions
$h_{\alpha \beta} (u) = (n_{\alpha}, n_{\beta})$, $1 \leq \alpha,
\beta \leq L$. For totally nonisotropic submanifolds we always have
$\det g_{ij} (u) \neq 0$ and $\det h_{\alpha \beta} (u) \neq 0$.
Note that usually in the local theory of submanifolds some
orthonormal bases in the normal spaces $\mathbb{N}_u$ are
considered, but it is fundamentally important for our approach to
consider arbitrary bases in the normal spaces $\mathbb{N}_u$.
Therefore, we develop such a general approach here and present in
detail the corresponding general fundamental relations, formulae,
and equations of the local theory of totally nonisotropic
submanifolds in pseudo-Euclidean spaces in the form necessary for
us.

\subsection{Gauss and Weingarten decompositions}

Since the set of vectors $(r_1 (u), \ldots, r_N (u), n_1 (u),
\ldots, n_L (u))$ forms a basis in $\mathbb{E}^{N + L}_q$ at each
point of the submanifold $M^N$, we can decompose each of the vectors
on the submanifold $M^N$ with respect to this basis, in particular,
the vectors $\pa^2 r / \pa u^i \pa u^j,$ $1 \leq i, j \leq N,$ and
the vectors $\pa n_{\alpha} / \pa u^i,$ $1 \leq \alpha \leq L,$ $1
\leq i \leq N$, getting the Gauss decomposition \be {\partial^2 r
\over
\partial u^i \partial u^j} = a^k_{ij} (u) {\partial r \over \partial
u^k} + b^{\beta}_{ij} (u) n_{\beta} (u) \label{1}\ee and the
Weingarten decomposition \be {\partial n_{\alpha} \over \partial
u^j} = c^k_{\alpha j} (u) {\partial r \over \partial u^k} +
d^{\beta}_{\alpha j} (u) n_{\beta} (u),\label{2} \ee where the
coefficients $a^k_{ij} (u),$ $b^{\beta}_{ij} (u),$ $c^k_{\alpha j}
(u)$, and $d^{\beta}_{\alpha j} (u)$ are smooth functions on the
submanifold $M^N$.

For each submanifold, there is a number of fundamental relations,
including the Gauss, Codazzi, and Ricci equations, between the
metric $g_{ij} (u)$, the functions $h_{\alpha \beta} (u)$, and the
coefficients $a^k_{ij} (u)$, $b^{\beta}_{ij} (u)$, $c^k_{\alpha j}
(u)$, and $d^{\beta}_{\alpha j} (u)$. If $g_{ij} (u)$, $h_{\alpha
\beta} (u)$, $a^k_{ij} (u)$, $b^{\beta}_{ij} (u)$, $c^k_{\alpha j}
(u)$, and $d^{\beta}_{\alpha j} (u)$ satisfy locally all these
relations, then by the Bonnet theorem there always exists a unique
(up to motion in the ambient pseudo-Euclidean space) submanifold
with these differential-geometric objects.

It follows immediately from the Gauss decomposition (\ref{1}) that
the coefficients $a^k_{ij} (u)$ and $b^{\beta}_{ij} (u)$ are
symmetric with respect to the lower indices: \be a^k_{ij} (u) =
a^k_{ji} (u), \label{3}\ee \be b^{\beta}_{ij} (u) = b^{\beta}_{ji}
(u). \label{4}\ee In addition, taking the scalar product of the
Gauss decomposition (\ref{1}) with the tangent vectors $r_l (u)$, $1
\leq l \leq N,$ we have \be \left ({\pa^2 r \over \pa u^i \pa u^j},
{\pa r \over \pa u^l} \right ) = a^k_{ij} (u) g_{kl} (u).
\label{5}\ee Differentiating the relation \be \left ({\pa r \over
\pa u^i}, {\pa r \over \pa u^j} \right ) = g_{ij} (u) \label{6}\ee
with respect to $u^s$, we get \be \left ({\pa^2 r \over \pa u^i \pa
u^s}, {\pa r \over \pa u^j} \right ) + \left ({\pa r \over \pa u^i},
{\pa^2 r \over \pa u^j \pa u^s} \right ) = {\pa g_{ij} \over \pa
u^s} \label{7}\ee or, taking into account (\ref{5}), \be a^k_{is}
(u) g_{kj} (u) + a^k_{js} (u) g_{ki} (u) = {\pa g_{ij} \over \pa
u^s}. \label{8}\ee In addition, rearranging the indices and taking
into account the symmetry of $a^k_{ij} (u)$ with respect to the
lower indices (\ref{3}), we have \be a^k_{ij} (u) g_{ks} (u) +
a^k_{js} (u) g_{ki} (u) = {\pa g_{is} \over \pa u^j}, \label{9}\ee
\be a^k_{is} (u) g_{kj} (u) + a^k_{ij} (u) g_{ks} (u) = {\pa g_{sj}
\over \pa u^i}. \label{10}\ee From (\ref{8}), (\ref{9}), and
(\ref{10}) we obtain \be a^k_{ij} (u) g_{ks} (u) = {1 \over 2} \left
({\pa g_{sj} \over \pa u^i} + {\pa g_{is} \over \pa u^j} - {\pa
g_{ij} \over \pa u^s} \right ), \label{11}\ee or \be a^k_{ij} (u) =
{1 \over 2} g^{ks} (u) \left ({\pa g_{sj} \over \pa u^i} + {\pa
g_{is} \over \pa u^j} - {\pa g_{ij} \over \pa u^s} \right ),
\label{11a}\ee where $g^{ks} (u)$ is the contravariant metric that
is the inverse of the first fundamental form $g_{ij} (u)$: $g^{ik}
(u) g_{kj} (u) = \delta^i_j;$ i.e., the coefficients $a^k_{ij} (u)$
are the coefficients of the Levi-Civita connection $\Gamma^k_{ij}
(u)$ of the metric $g_{ij} (u)$: \be \Gamma^k_{ij} (u) = {1 \over 2}
g^{ks} (u) \left ({\pa g_{sj} \over \pa u^i} + {\pa g_{is} \over \pa
u^j} - {\pa g_{ij} \over \pa u^s} \right ) = a^k_{ij} (u).
\label{11b} \ee

Taking the scalar product of the Weingarten decomposition (\ref{2})
with the normal vectors $n_{\beta} (u)$, $1 \leq \beta \leq L,$ we
have \be \left ({\pa n_{\alpha} \over \pa u^i }, n_{\beta} \right )
= d^{\gamma}_{\alpha i} (u) h_{\gamma \beta} (u). \label{12}\ee The
expressions \be \varkappa_{\alpha \beta i} (u)  = \left ({\pa
n_{\alpha} \over \pa u^i }, n_{\beta} \right ) = d^{\gamma}_{\alpha
i} (u) h_{\gamma \beta} (u), \ \ \ \ 1 \leq \alpha, \beta \leq L,
\label{12a}\ee for any $\alpha$ and $\beta$ give components of a
covariant vector field and are called the {\it torsion coefficients
of the submanifold\/}, and the 1-forms $\varkappa_{\alpha \beta i}
(u) d u^i$, $1 \leq \alpha, \beta \leq L$, are called the {\it
torsion forms of the submanifold}.

Differentiating the relation \be (n_{\alpha}, n_{\beta})  =
h_{\alpha \beta} (u) \label{13}\ee with respect to $u^i$, we get \be
\left ({\pa n_{\alpha} \over \pa u^i}, n_{\beta} \right ) + \left
(n_{\alpha}, {\pa n_{\beta} \over \pa u^i} \right ) = {\pa h_{\alpha
\beta} \over \pa u^i}, \label{14}\ee and taking into account
(\ref{12}), \be d^{\gamma}_{\alpha i} (u) h_{\gamma \beta} (u) +
d^{\gamma}_{\beta i} (u) h_{\gamma \alpha} (u) = {\pa h_{\alpha
\beta} \over \pa u^i}, \label{15}\ee or \be \varkappa_{\alpha \beta
i} (u) + \varkappa_{\beta \alpha i} (u) = {\pa h_{\alpha \beta}
\over \pa u^i}. \label{15a}\ee

Taking the scalar product of the Gauss decomposition (\ref{1}) with
the normal vectors $n_{\alpha} (u)$, $1 \leq \alpha \leq L,$ we have
\be \left ({\pa^2 r \over \pa u^i \pa u^j}, n_{\alpha} (u) \right )
= b^{\beta}_{ij} (u) h_{\beta \alpha} (u), \label{16}\ee and taking
the scalar product of the Weingarten decomposition (\ref{2}) with
the tangent vectors $r_j (u)$, $1 \leq j \leq N,$ we get \be \left
({\pa n_{\alpha} \over \pa u^i }, {\pa r \over \pa u^j} \right ) =
c^k_{\alpha i} (u) g_{kj} (u). \label{17}\ee

The expressions \be \omega_{\alpha ij} (u) = \left ({\pa^2 r \over
\pa u^i \pa u^j}, n_{\alpha} (u) \right ) = b^{\beta}_{ij} (u)
h_{\beta \alpha} (u), \ \ \ \ 1 \leq \alpha \leq L, \label{kvf} \ee
for any $\alpha$ give components of a symmetric covariant tensor and
are called the {\it second fundamental forms of the submanifold\/}.
To each basis vector $n_{\alpha} (u)$, $1 \leq \alpha \leq L,$ of
the normal space of the submanifold there corresponds its own second
fundamental form.

Differentiating the relation \be \left ({\pa r \over \pa u^i},
n_{\alpha} (u) \right )  = 0 \label{18}\ee with respect to $u^j$, we
get \be \left ({\pa^2 r \over \pa u^i \pa u^j}, n_{\alpha} (u)
\right ) + \left ({\pa r \over \pa u^i}, {\pa n_{\alpha} \over \pa
u^j} \right )= 0, \label{19}\ee or, taking into account (\ref{16})
and (\ref{17}), \be b^{\beta}_{ij} (u) h_{\beta \alpha} (u) +
c^k_{\alpha j} (u) g_{ki} (u)= 0. \label{20}\ee

Thus, the following relation always holds: \be c^k_{\alpha i} (u) =
- g^{ks} (u) b^{\beta}_{si} (u) h_{\beta \alpha} (u) = - g^{ks} (u)
\omega_{\alpha si} (u). \label{21}\ee For any $\alpha$ the
coefficients  $c^k_{\alpha i} (u)$, $1 \leq \alpha \leq L$, give
components of a tensor of rank (1, 1) (an affinor). These affinors
are called the {\it Weingarten operators of the submanifold\/}.

\subsection{Consistency conditions}

The consistency conditions for the Gauss decomposition (\ref{1}) and
the Weingarten decomposition (\ref{2}) are represented by
fundamental equations in submanifold theory, namely, by the Gauss
equations, the Codazzi equations, and the Ricci equations. Here, we
consider the consistency conditions for the Gauss decomposition
(\ref{1}) and the Weingarten decomposition (\ref{2}) in the form
necessary for us.

Differentiating the Gauss decomposition (\ref{1}) with respect to
$u^s$, we find \be {\partial^3 r \over \partial u^i
\partial u^j \partial u^s} = {\pa {\Gamma}^k_{ij} \over \pa u^s}
{\partial r \over
\partial u^k} + {\Gamma}^k_{ij} (u)
{\partial^2 r \over
\partial u^k \partial u^s} + {\pa b^{\beta}_{ij} \over \pa u^s} n_{\beta} (u) +
b^{\beta}_{ij} (u) {\pa n_{\beta} \over \pa u^s}. \label{22}\ee
Using the Gauss decomposition (\ref{1}) and the Weingarten
decomposition (\ref{2}), we obtain \bea &&{\partial^3 r \over
\partial u^i \partial u^j
\partial u^s} = {\pa {\Gamma}^k_{ij} \over \pa u^s} {\partial r
\over
\partial u^k} + {\Gamma}^k_{ij} (u) \left ( {\Gamma}^l_{ks} (u)
{\partial r \over \partial u^l} + b^{\beta}_{ks} (u) n_{\beta} (u)
\right ) + \nn\\ && + {\pa b^{\beta}_{ij} \over \pa u^s} n_{\beta}
(u) + b^{\beta}_{ij} (u) \left (c^k_{\beta s} (u) {\partial r \over
\partial u^k} + d^{\gamma}_{\beta s} (u) n_{\gamma} (u)\right ).
\label{23}\eea The condition of symmetry with respect to the indices
$j$ and $s$ yields \bea &&{\pa {\Gamma}^k_{ij} \over \pa u^s}
{\partial r \over
\partial u^k} + {\Gamma}^l_{ij} (u) \left ( {\Gamma}^k_{ls} (u)
{\partial r \over \partial u^k} + b^{\beta}_{ls} (u) n_{\beta} (u)
\right ) + \nn\\ && + {\pa b^{\beta}_{ij} \over \pa u^s} n_{\beta}
(u) + b^{\gamma}_{ij} (u) \left (c^k_{\gamma s} (u) {\partial r
\over
\partial u^k} + d^{\beta}_{\gamma s} (u) n_{\beta} (u)\right ) = \nn\\
&& = {\pa {\Gamma}^k_{is} \over \pa u^j} {\partial r \over
\partial u^k} + {\Gamma}^l_{is} (u) \left ( {\Gamma}^k_{lj} (u)
{\partial r \over \partial u^k} + b^{\beta}_{lj} (u) n_{\beta} (u)
\right ) + \nn\\ && + {\pa b^{\beta}_{is} \over \pa u^j} n_{\beta}
(u) + b^{\gamma}_{is} (u) \left (c^k_{\gamma j} (u) {\partial r
\over
\partial u^k} + d^{\beta}_{\gamma j} (u) n_{\beta} (u)\right ). \label{24}\eea

The coefficients of $\pa r/\pa u^k$ give the Gauss equations \be
{\pa {\Gamma}^k_{ij} \over \pa u^s} - {\pa {\Gamma}^k_{is} \over \pa
u^j} + {\Gamma}^l_{ij} (u) {\Gamma}^k_{ls} (u) - {\Gamma}^l_{is} (u)
{\Gamma}^k_{lj} (u) =  b^{\gamma}_{is} (u) c^k_{\gamma j} (u) -
b^{\gamma}_{ij} (u) c^k_{\gamma s} (u), \label{25}\ee or \be
R^k_{isj} (u) = b^{\gamma}_{is} (u) c^k_{\gamma j} (u) -
b^{\gamma}_{ij} (u) c^k_{\gamma s} (u), \label{26}\ee where
$R^k_{isj} (u)$ is the Riemannian curvature tensor of the metric
$g_{ij} (u)$.

The coefficients of $n_{\beta} (u)$ give the Codazzi equations \be
{\Gamma}^l_{ij} (u)  b^{\beta}_{ls} (u) + {\pa b^{\beta}_{ij} \over
\pa u^s} + b^{\gamma}_{ij} (u) d^{\beta}_{\gamma s} (u) =
{\Gamma}^l_{is} (u) b^{\beta}_{lj} (u) + {\pa b^{\beta}_{is} \over
\pa u^j} + b^{\gamma}_{is} (u) d^{\beta}_{\gamma j} (u).
\label{27}\ee

Differentiating the Weingarten decomposition (\ref{2}) with respect
to $u^s$, we have \be {\partial^2 n_{\alpha} \over \partial u^j
\partial u^s} = {\pa c^k_{\alpha j} \over \pa u^s} {\partial r \over
\partial u^k} + c^k_{\alpha j} (u) {\partial^2 r \over
\partial u^k \partial u^s} + {\pa d^{\beta}_{\alpha j} \over \pa
u^s} n_{\beta} (u) + d^{\beta}_{\alpha j} (u) {\pa n_{\beta} \over
\pa u^s}. \label{28} \ee Using the Gauss decomposition (\ref{1}) and
the Weingarten decomposition (\ref{2}), we get \bea &&{\partial^2
n_{\alpha} \over
\partial u^j
\partial u^s} = {\pa c^k_{\alpha j} \over \pa u^s} {\partial r \over
\partial u^k} + c^k_{\alpha j} (u) \left ( {\Gamma}^l_{ks} (u)
{\partial r \over
\partial u^l} + b^{\beta}_{ks} (u) n_{\beta} (u) \right ) + \nn\\
&& + {\pa d^{\beta}_{\alpha j} \over \pa u^s} n_{\beta} (u) +
d^{\beta}_{\alpha j} (u) \left (c^k_{\beta s} (u) {\partial r \over
\partial u^k} + d^{\gamma}_{\beta s} (u) n_{\gamma} (u)\right ). \label{29}
\eea The condition of symmetry with respect to the indices $j$ and
$s$ gives \bea &&{\pa c^k_{\alpha j} \over \pa u^s} {\partial r
\over
\partial u^k} + c^l_{\alpha j} (u) \left ( {\Gamma}^k_{ls} (u)
{\partial r \over
\partial u^k} + b^{\beta}_{ls} (u) n_{\beta} (u) \right ) + \nn\\
&& + {\pa d^{\beta}_{\alpha j} \over \pa u^s} n_{\beta} (u) +
d^{\gamma}_{\alpha j} (u) \left (c^k_{\gamma s} (u) {\partial r
\over
\partial u^k} + d^{\beta}_{\gamma s} (u) n_{\beta} (u)\right ) = \nn\\
&& = {\pa c^k_{\alpha s} \over \pa u^j} {\partial r \over
\partial u^k} + c^l_{\alpha s} (u) \left ( {\Gamma}^k_{lj} (u)
{\partial r \over
\partial u^k} + b^{\beta}_{lj} (u) n_{\beta} (u) \right ) + \nn\\
&& + {\pa d^{\beta}_{\alpha s} \over \pa u^j} n_{\beta} (u) +
d^{\gamma}_{\alpha s} (u) \left (c^k_{\gamma j} (u) {\partial r
\over
\partial u^k} + d^{\beta}_{\gamma j} (u) n_{\beta} (u)\right ). \label{30} \eea
The coefficients of $n_{\beta} (u)$ give the Ricci equations \be
c^l_{\alpha j} (u) b^{\beta}_{ls} (u) + {\pa d^{\beta}_{\alpha j}
\over \pa u^s} + d^{\gamma}_{\alpha j} (u) d^{\beta}_{\gamma s} (u)
 = c^l_{\alpha s} (u) b^{\beta}_{lj} (u) +
 {\pa d^{\beta}_{\alpha s} \over \pa u^j} +
d^{\gamma}_{\alpha s} (u) d^{\beta}_{\gamma j} (u). \label{31} \ee
The coefficients of $\pa r/\pa u^k$ give the Codazzi equations \be
{\pa c^k_{\alpha j} \over \pa u^s} + c^l_{\alpha j} (u)
{\Gamma}^k_{ls} (u) + d^{\gamma}_{\alpha j} (u) c^k_{\gamma s} (u) =
{\pa c^k_{\alpha s} \over \pa u^j} + c^l_{\alpha s} (u)
{\Gamma}^k_{lj} (u) + d^{\gamma}_{\alpha s} (u) c^k_{\gamma j} (u).
\label{32} \ee

If relations (\ref{21}), (\ref{3}), (\ref{8}), (\ref{11b}), and
(\ref{15}) hold, then equations (\ref{32}) are equivalent to the
Codazzi equations (\ref{27}). Indeed, substituting $c^k_{\alpha i}
(u)$ from (\ref{21}) into (\ref{32}), we obtain \bea && - {\pa
g^{ks} \over \pa u^j} b^{\beta}_{si} (u) h_{\beta \alpha} (u) -
g^{ks} (u) {\pa b^{\beta}_{si} \over \pa u^j} h_{\beta \alpha} (u) -
g^{ks} (u)
b^{\beta}_{si} (u) {\pa h_{\beta \alpha} \over \pa u^j} - \nn\\
&& - g^{ls} (u) b^{\beta}_{si} (u) h_{\beta \alpha} (u)
{\Gamma}^k_{lj} (u) - d^{\beta}_{\alpha i} (u) g^{ks} (u)
b^{\gamma}_{sj} (u) h_{\gamma \beta} (u) = \nn\\ && = - {\pa g^{ks}
\over \pa u^i} b^{\beta}_{sj} (u) h_{\beta \alpha} (u) - g^{ks} (u)
{\pa b^{\beta}_{sj} \over \pa u^i} h_{\beta \alpha} (u) - g^{ks} (u)
b^{\beta}_{sj} (u) {\pa h_{\beta \alpha} \over \pa u^i} - \nn\\ && -
g^{ls} (u) b^{\beta}_{sj} (u) h_{\beta \alpha} (u) {\Gamma}^k_{li}
(u) - d^{\beta}_{\alpha j} (u) g^{ks} (u) b^{\gamma}_{si} (u)
h_{\gamma \beta} (u). \label{33} \eea From the compatibility
condition of the connection $\Gamma^r_{pj} (u)$ with the metric
$g_{rp} (u)$ (or from relations (\ref{3}), (\ref{8}), and
(\ref{11b})), for the derivative of the contravariant metric $g^{ks}
(u)$ we have  \bea && {\pa g^{ks} \over \pa u^j} = - g^{kp} (u)
\left ( \Gamma^r_{pj} (u) g_{rl} (u) + \Gamma^r_{lj} (u) g_{rp} (u)
\right ) g^{ls} (u) = \nn\\ && = - \Gamma^s_{pj} (u) g^{kp} (u)-
\Gamma^k_{lj} (u) g^{ls} (u). \label{34} \eea

Using relations (\ref{34}) and (\ref{15}), from (\ref{33}) we get
\bea && g^{kp} (u) \left ( \Gamma^r_{pj} (u) g_{rl} (u) +
\Gamma^r_{lj} (u) g_{rp} (u) \right ) g^{ls} (u) b^{\beta}_{si} (u)
h_{\beta \alpha} (u) - \nn\\ && - g^{ks} (u) {\pa b^{\beta}_{si}
\over \pa u^j} h_{\beta \alpha} (u) - g^{ks} (u) b^{\beta}_{si} (u)
\left ( d^{\gamma}_{\alpha j} (u) h_{\gamma
\beta} (u) + d^{\gamma}_{\beta j} (u) h_{\gamma \alpha} (u) \right ) - \nn\\
&& - g^{ls} (u) b^{\beta}_{si} (u) h_{\beta \alpha} (u)
{\Gamma}^k_{lj} (u) - d^{\beta}_{\alpha i} (u) g^{ks} (u)
b^{\gamma}_{sj} (u) h_{\gamma \beta} (u) = \nn\\ && = g^{kp} (u)
\left ( \Gamma^r_{pi} (u) g_{rl} (u) + \Gamma^r_{li} (u) g_{rp} (u)
\right ) g^{ls} (u)b^{\beta}_{sj} (u) h_{\beta \alpha} (u) - \nn\\
&& - g^{ks} (u) {\pa b^{\beta}_{sj} \over \pa u^i} h_{\beta \alpha}
(u) - g^{ks} (u) b^{\beta}_{sj} (u) \left ( d^{\gamma}_{\alpha i}
(u) h_{\gamma \beta} (u) + d^{\gamma}_{\beta i} (u) h_{\gamma
\alpha} (u) \right ) - \nn\\ && - g^{ls} (u) b^{\beta}_{sj} (u)
h_{\beta \alpha} (u) {\Gamma}^k_{li} (u) - d^{\beta}_{\alpha j} (u)
g^{ks} (u) b^{\gamma}_{si} (u) h_{\gamma \beta} (u), \label{35} \eea
or \be \Gamma^s_{pj} (u) b^{\beta}_{si} (u) - {\pa b^{\beta}_{pi}
\over \pa u^j} - b^{\gamma}_{pi} (u) d^{\beta}_{\gamma j} (u)
 = \Gamma^s_{pi} (u) b^{\beta}_{sj} (u) - {\pa b^{\beta}_{pj} \over
\pa u^i} - b^{\gamma}_{pj} (u) d^{\beta}_{\gamma i} (u), \label{36}
\ee i.e., the Codazzi equations (\ref{27}).

\subsection{Bonnet theorem}

For totally nonisotropic submanifolds in pseudo-Euclidean spaces, an
analog of the classical Bonnet theorem holds. Let a
pseudo-Riemannian metric $g_{ij} (u) du^i du^j$, symmetric 2-forms
$\omega_{\alpha ij} (u) du^i du^j$, $\omega_{\alpha ij} (u) =
\omega_{\alpha ji} (u)$, $1 \leq \alpha \leq L,$ 1-forms
$\varkappa_{\alpha \beta i} (u) du^i$, $1 \leq \alpha, \beta \leq
L,$ and functions $h_{\alpha \beta} (u)$, $1 \leq \alpha, \beta \leq
L$, such that $\det h_{\alpha \beta} (u) \neq 0$ and $h_{\alpha
\beta} (u) = h_{\beta \alpha} (u)$, $1 \leq \alpha, \beta \leq L$,
be locally given. If in this case relations (\ref{15a}) as well as
the Gauss equations {\rm (\ref{26})}, the Codazzi equations {\rm
(\ref{27})} and the Ricci equations {\rm (\ref{31})} are satisfied
for the forms $g_{ij} (u)$, $\omega_{\alpha ij} (u)$,
$\varkappa_{\alpha \beta i} (u)$ and the functions $h_{\alpha \beta}
(u)$ (the coefficients $b^{\beta}_{ij} (u)$, $c^k_{\alpha i} (u)$,
and $d^{\beta}_{\alpha j} (u)$ are uniquely determined by formulae
(\ref{kvf}), (\ref{21}), and (\ref{12a}), respectively), then there
exists a unique {\rm (}up to motions in the ambient pseudo-Euclidean
space{\rm )} smooth totally nonisotropic $N$-dimensional submanifold
$M^N$ with the first fundamental form $d s^2 = g_{ij} (u) du^i
du^j$, the Gram matrix $h_{\alpha \beta} (u)$, $1 \leq \alpha, \beta
\leq L$, of scalar products of the basis vectors in the normal
spaces, the second fundamental forms $\omega_{\alpha ij} (u) d u^i d
u^j$, and the torsion forms $\varkappa_{\alpha \beta i} (u) d u^i$
in an $(N + L)$-dimensional pseudo-Euclidean space, the signature of
which is determined by the sum of the signatures of the metrics
$g_{ij} (u)$, $1 \leq i, j \leq N$, and $h_{\alpha \beta} (u)$, $1
\leq \alpha, \beta \leq L$.

\section{Submanifolds with zero torsion in \\ pseudo-Eucli\-dean spaces}

Let us consider the class of totally nonisotropic smooth
$N$-dimensional submanifolds with {\it zero tor\-sion\/} in $(N +
L)$-dimensional pseudo-Euclidean spaces; i.e., all the torsion forms
$d^{\beta}_{\alpha i} (u) d u^i$, $1 \leq \alpha, \beta \leq L$, of
submanifolds of this class vanish, $d^{\beta}_{\alpha i} (u) = 0$,
$1 \leq \alpha, \beta \leq L$, $1 \leq i \leq N$, for the chosen
bases in the normal spaces. In this case it follows immediately from
relations (\ref{15}) that the functions $h_{\alpha \beta} (u)$, $1
\leq \alpha, \beta \leq L$, must be constant: $h_{\alpha \beta} (u)
= \mu_{\alpha \beta}$, $\mu_{\alpha \beta} = {\rm \ const}$, where
$\mu_{\alpha \beta} = \mu_{\beta \alpha}$ and $\det (\mu_{\alpha
\beta}) \neq 0$ by virtue of the definition of these functions.
Relations (\ref{15}) hold in this case. Note that if the functions
$h_{\alpha \beta} (u)$, $1 \leq \alpha, \beta \leq L$, are constant,
then relations (\ref{15}) are equivalent to the condition \be
d^{\gamma}_{\alpha i} (u) h_{\gamma \beta} (u) + d^{\gamma}_{\beta
i} (u) h_{\gamma \alpha} (u) = 0, \label{1z}\ee i.e., the
skew-symmetry condition of the torsion 1-forms $\varkappa_{\alpha
\beta i} (u) d u^i = d^{\gamma}_{\alpha i} (u) h_{\gamma \beta} (u)
d u^i$ with respect to the indices $\alpha$ and $\beta$:
$\varkappa_{\alpha \beta i} (u) = - \varkappa_{\beta \alpha i} (u)$.
The converse is also true; i.e., the functions $h_{\alpha \beta}
(u)$, $1 \leq \alpha, \beta \leq L$, are constant if and only if the
torsion 1-forms $\varkappa_{\alpha \beta i} (u) d u^i$ are
skew-symmetric with respect to the indices $\alpha$ and $\beta$.

For submanifolds with zero torsion, the following relations hold:
\be c^k_{\alpha i} (u) = - g^{ks} (u) b^{\beta}_{si} (u) \mu_{\beta
\alpha}, \ \ \ \ \omega_{\alpha ij} (u) = b^{\beta}_{ij} (u)
\mu_{\beta \alpha}, \label{2z}\ee the Gauss equations \be R^k_{isj}
(u) = b^{\gamma}_{is} (u) c^k_{\gamma j} (u) - b^{\gamma}_{ij} (u)
c^k_{\gamma s} (u), \label{3z}\ee the Codazzi equations \be
{\Gamma}^l_{ij} (u) b^{\beta}_{lk} (u) + {\pa b^{\beta}_{ij} \over
\pa u^k} = {\Gamma}^l_{ik} (u) b^{\beta}_{lj} (u) + {\pa
b^{\beta}_{ik} \over \pa u^j}, \label{4z}\ee and the Ricci equations
\be c^l_{\alpha j} (u) b^{\beta}_{ls} (u)
 = c^l_{\alpha s} (u) b^{\beta}_{lj} (u). \label{5z} \ee

The Codazzi equations (\ref{4z}) can be rewritten in the form \be
\nabla_k b^{\alpha}_{ij} = \nabla_j b^{\alpha}_{ik} \label{6z} \ee
or \be \nabla_k \omega_{\alpha ij} = \nabla_j \omega_{\alpha ik},
\label{6za} \ee where $\nabla$ is the covariant derivative generated
by the Levi-Civita connection of the first fundamental form $g_{ij}
(u)$.

Using relation (\ref{2z}), one can rewrite the Gauss equations
(\ref{3z}) in the form \be R_{ijkl} (u) = b^{\alpha}_{ik} (u)
\mu_{\alpha \beta} b^{\beta}_{lj} (u) - b^{\alpha}_{jk} (u)
\mu_{\alpha \beta} b^{\beta}_{li} (u) \label{7z} \ee or \be
R^{ik}_{sj} (u) = c^i_{\alpha j} (u) \mu^{\alpha \beta} c^k_{\beta
s} (u) - c^i_{\alpha s} (u) \mu^{\alpha \beta} c^k_{\beta j} (u),
\label{7zab} \ee and also \be R_{ijkl} (u) = \omega_{\alpha ik} (u)
\mu^{\alpha \beta} \omega_{\beta lj} (u) - \omega_{\alpha jk} (u)
\mu^{\alpha \beta} \omega_{\beta li} (u), \label{7za} \ee where the
matrix $\mu^{\alpha \beta}$ is the inverse of the matrix
$\mu_{\alpha \beta}$: $ \mu^{\alpha \gamma} \mu_{\gamma \beta}=
\delta^{\alpha}_{\beta}$, and the Ricci equations (\ref{5z}) take
the form \be b^{\alpha}_{ik} (u) g^{kl} (u) b^{\beta}_{lj} (u) -
b^{\alpha}_{jk} (u) g^{kl} (u) b^{\beta}_{li} (u) = 0 \label{8z} \ee
or \be c^i_{\alpha j} (u) g_{ik} (u) c^k_{\beta s} (u) - c^i_{\alpha
s} (u) g_{ik} (u) c^k_{\beta j} (u) = 0, \label{8zab} \ee and also
\be \omega_{\alpha ik} (u) g^{kl} (u) \omega_{\beta lj} (u) -
\omega_{\alpha jk} (u) g^{kl} (u) \omega_{\beta li} (u) = 0.
\label{8za} \ee

\subsection{Flat submanifolds with zero torsion}

Now we consider flat submanifolds with zero torsion in
pseudo-Euclidean spaces, i.e., torsionless submanifolds with flat
metrics, namely, with flat first fundamental forms $g_{ij} (u)$ on
the submanifolds. In this case, we can consider that $u = (u^1,
\ldots, u^N)$ are some flat coordinates of the metric $g_{ij} (u)$.
In flat coordinates the metric is a constant nondegenerate symmetric
matrix $\eta_{ij}$, $\eta_{ij} = {\rm \ const}$, $\det \,
(\eta_{ij}) \neq 0$, $\eta_{ij} = \eta_{ji}$, and the Codazzi
equations (\ref{4z}), (\ref{6z}) take the form \be {\pa
b^{\alpha}_{ij} \over \pa u^k} = {\pa b^{\alpha}_{ik} \over \pa u^j}
\label{9z} \ee or \be {\pa \omega_{\alpha ij} \over \pa u^k} = {\pa
\omega_{\alpha ik} \over \pa u^j}. \label{9za} \ee Hence, there
locally exist some functions $\psi_{\alpha} (u)$, $1 \leq \alpha
\leq L$, such that \be \omega_{\alpha ij} (u) = {\pa^2
{\psi}_{\alpha} \over \pa u^i \pa u^j}. \label{10z} \ee We have thus
proved the following important lemma. \bl {\rm \cite{3m}}
\label{lem} All the second fundamental forms of each flat
torsionless submani\-fold in a pseudo-Euclidean space are locally
Hessians in any flat coordinates on the submanifold. \el

It follows from Lemma \ref{lem} that in any flat coordinates the
Gauss equations (\ref{7za}) have the form \be \sum_{\alpha = 1}^L
\sum_{\beta = 1}^L \mu^{\alpha \beta} \left ({\pa^2 \psi_{\alpha}
\over \pa u^i \pa u^k} {\pa^2 \psi_{\beta} \over \pa u^j \pa u^l} -
{\pa^2 \psi_{\alpha} \over \pa u^i \pa u^l} {\pa^2 \psi_{\beta}
\over \pa u^j \pa u^k} \right ) = 0, \label{7zb} \ee and the Ricci
equations (\ref{8za}) have the form \be \sum_{i = 1}^N \sum_{j =
1}^N \eta^{ij} \left ({\pa^2 \psi_{\alpha} \over \pa u^i \pa u^k}
{\pa^2 \psi_{\beta} \over \pa u^j \pa u^l} - {\pa^2 \psi_{\alpha}
\over \pa u^i \pa u^l} {\pa^2 \psi_{\beta} \over \pa u^j \pa u^k}
\right ) = 0,  \label{8zb} \ee where the matrix $\eta^{ij}$ is the
inverse of the matrix $\eta_{ij}$: $\eta^{is} \eta_{sj} =
\delta^i_j$.

\bt {\rm \cite{1}, \cite{3m}, \cite{2m}} The class of
$N$-dimensional flat torsionless submanifolds in $(N +
L)$-dimensional pseudo-Euclidean spaces is described {\rm (}in flat
coordinates{\rm )} by the system of nonlinear equations {\rm
(\ref{7zb}), (\ref{8zb})} for functions $\psi_{\alpha} (u),$ $1 \leq
\alpha \leq L.$ Here, $\eta^{ij}$, $1 \leq i, j \leq N$, and
$\mu^{\alpha \beta}$, $1 \leq \alpha, \beta \leq L$, are arbitrary
constant nondegenerate symmetric matrices, $\eta^{ij} = \eta^{ji},$
$\eta^{ij} = {\rm const},$ $\det (\eta^{ij}) \neq 0$, $\mu^{\alpha
\beta}= {\rm const}$, $\mu^{\alpha \beta} = \mu^{\beta \alpha}$,
$\det (\mu^{\alpha \beta}) \neq 0${\rm ;} the signature of the
ambient $(N + L)$-dimensional pseudo-Euclidean space is the sum of
the signatures of the metrics $\eta^{ij}$ and $\mu^{\alpha
\beta}${\rm ;} ${\bf I} = d s^2 = \eta_{ij} d u^i d u^j$ is the
first fundamental form, where $\eta_{ij}$ is the inverse of the
matrix $\eta^{ij}$, $\eta^{is} \eta_{sj} = \delta^i_j$, and ${\bf
II}_{\alpha} = (\pa^2 \psi_{\alpha} / (\pa u^i \pa u^j)) d u^i d
u^j,$ $1 \leq \alpha \leq L,$ are the second fundamental forms given
by the Hessians of the functions $\psi_{\alpha} (u),$ $1 \leq \alpha
\leq L$, for the corresponding flat torsionless submanifold
determined by an arbitrary solution of the system of nonlinear
equations {\rm (\ref{7zb}), (\ref{8zb})}. \et

According to the Bonnet theorem, any solution $\psi_{\alpha} (u),$
$1 \leq \alpha \leq L,$ of the nonlinear system (\ref{7zb}),
(\ref{8zb}) determines a unique (up to motion in the ambient
pseudo-Euclidean space) totally nonisotropic $N$-dimensional flat
torsionless submanifold in the cor\-responding $(N + L)$-dimensional
pseudo-Euclidean space with the first fundamental form $\eta_{ij} d
u^i d u^j$ and the second fundamental forms $\omega_{\alpha} (u) =
(\pa^2 \psi_{\alpha} / (\pa u^i \pa u^j)) d u^i d u^j$, $1 \leq
\alpha \leq L$, given by the Hessians of the functions
$\psi_{\alpha} (u),$ $1 \leq \alpha \leq L$, and the constant Gram
matrix $\mu_{\alpha \beta}$, $1 \leq \alpha, \beta \leq L$, of the
scalar products of basic vectors in the normal spaces. It is obvious
that we can always add arbitrary terms linear in the coordinates
$u^1, \ldots, u^N$ to any solution of the system (\ref{7zb}),
(\ref{8zb}), but the set of the second fundamental forms and the
corresponding submanifold will remain the same. Moreover, any two
sets of the second fundamental forms of the form $\omega_{\alpha ij}
(u) = \pa^2 \psi_{\alpha} / (\pa u^i \pa u^j)$, $1 \leq \alpha \leq
L$, coincide if and only if the corresponding functions
$\psi_{\alpha} (u),$ $1 \leq \alpha \leq L$, coincide up to terms
linear in the coordinates; hence we must not distinguish here
solutions of the nonlinear system (\ref{7zb}), (\ref{8zb}) that
differ by terms linear in the coordinates $u^1, \ldots, u^N$.

We consider the following linear problem with parameters for vector
functions $\pa a(u) / \pa u^i$, $1 \leq i \leq N,$ and $b_{\alpha}
(u)$, $1 \leq \alpha \leq L$: \be {\pa^2 a \over \pa u^i \pa u^j} =
\lambda \, \mu^{\alpha \beta} \omega_{\alpha ij} (u) b_{\beta} (u),\
\ \ \ {\pa b_{\alpha} \over \pa u^i} = \rho \, \eta^{kj}
\omega_{\alpha ij} (u) {\pa a \over \pa u^k}, \label{t1} \ee where
$\eta^{ij}$, $ 1 \leq i, j \leq N,$ and $\mu^{\alpha \beta},$ $1
\leq \alpha, \beta \leq L$, are arbitrary constant nondegenerate
symmetric matrices, $\eta^{ij} = \eta^{ji},$ $\eta^{ij} = {\rm
const},$ $\det (\eta^{ij}) \neq 0$, $\mu^{\alpha \beta}= {\rm
const}$, $\mu^{\alpha \beta} = \mu^{\beta \alpha}$, $\det
(\mu^{\alpha \beta}) \neq 0$; $\lambda$ and $\rho$ are arbitrary
constants (parameters) \cite{1}. (In fact, only one of the
parameters $\lambda$ and $\rho$ is essential.) It is obvious that
the coefficients $\omega_{\alpha ij} (u)$, $1 \leq \alpha \leq L,$
here must be symmetric matrix functions, $\omega_{\alpha ij} (u) =
\omega_{\alpha ji} (u)$.

The consistency conditions for the linear system (\ref{t1}) are
equivalent to the nonli\-near system (\ref{7zb}), (\ref{8zb})
describing the class of $N$-dimensional flat torsionless
subma\-nifolds in $(N + L)$-dimensional pseudo-Euclidean spaces.

\bt  {\rm \cite{1}} The nonlinear system {\rm (\ref{7zb}),
(\ref{8zb})} is integrable by the inverse scattering method. \et

\bd {\rm A class of submanifolds in a Euclidean or pseudo-Euclidean
space is called {\it integrable} if the system of the fundamental
Gauss--Codazzi--Ricci equations giving this class of submanifolds is
integrable.}\ed

Essentially, the theory of integrable classes of surfaces in
${\mathbb{E}^3}$ goes back to the classical differential geometry of
the XIX century, when there were established remarkable properties
of some nonlinear partial differential equations (in particular, the
sine-Gordon equation and the Liouville equation) describing some
important classes of surfaces in ${\mathbb{E}^3}$. From the modern
viewpoint, after the methods of the soliton theory were discovered
and worked out and the theory of integrable nonlinear partial
differential equations was developed in the second half of the XX
century, it is clear that these properties indicate the
integrability of these nonlinear equations by the inverse scattering
method. In connection with the rapid and intensive development of
the theory of integrable systems, integrable classes of surfaces
have been considered and studied in many papers; in particular, we
mention the cycle of Sym's papers (see \cite{sym1}, \cite{sym2}) and
also the papers \cite{sav}--\cite{fok}. We also note that the
considered notion of integrability concerns only classes of surfaces
or submanifolds and makes no sense for a single surface or
submanifold. In particular, the definition of an {\it integrable
surface} via the integrability of {\it its} Gauss--Codazzi--Ricci
equations in \cite{bob} and \cite{fok} is quite absurd, since for
any surface {\it its} Gauss--Codazzi--Ricci equations are always
satisfied identically. In the context of the integrability of the
Gauss--Codazzi--Ricci equations, one can only speak of integrable
classes of surfaces and of whether a given surface belongs to a
certain integrable class of surfaces, but not of integrable
Gauss--Codazzi--Ricci equations of a concrete surface. Of course, a
concrete surface (or a submanifold) can belong to various integrable
and nonintegrable classes.

\bt The class of flat torsionless submanifolds in any Euclidean or
pseudo-Euclidean space is integrable. \et

\section{$k$-potential reductions and $k$-potential \\ submani\-folds
in pseudo-Euclidean spaces}

Consider the case when $L = k N + p$, where $k$ is an arbitrary
positive integer and $p$ is an arbitrary nonnegative integer, $p
\geq 0$.

In this case, the Gauss equations (\ref{7zb}) and the Ricci
equations (\ref{8zb}) can be rewritten in the following form: \bea
&& \sum_{r = 1}^k \sum_{s = 1}^k \sum_{\alpha = (r - 1) N + 1}^{r N}
\sum_{\beta = (s - 1) N + 1}^{s N} \mu^{\alpha \beta} \left ({\pa^2
\psi_{\alpha} \over \pa u^i \pa u^q} {\pa^2 \psi_{\beta} \over \pa
u^j \pa u^l} - {\pa^2 \psi_{\alpha} \over \pa u^i \pa u^l} {\pa^2
\psi_{\beta} \over \pa
u^j \pa u^q} \right ) + \nn\\
&& + \sum_{\alpha = k N + 1}^L \sum_{\beta = 1}^L \mu^{\alpha \beta}
\left ({\pa^2 \psi_{\alpha} \over \pa u^i \pa u^q} {\pa^2
\psi_{\beta} \over \pa u^j \pa u^l} - {\pa^2 \psi_{\alpha} \over \pa
u^i \pa u^l} {\pa^2 \psi_{\beta} \over \pa u^j \pa u^q} \right ) +
\nn\\ && +  \sum_{\alpha = 1}^{k N} \sum_{\beta = k N + 1}^L
\mu^{\alpha \beta} \left ({\pa^2 \psi_{\alpha} \over \pa u^i \pa
u^q} {\pa^2 \psi_{\beta} \over \pa u^j \pa u^l} - {\pa^2
\psi_{\alpha} \over \pa u^i \pa u^l} {\pa^2 \psi_{\beta} \over \pa
u^j \pa u^q} \right ) = 0, \label{g3a} \eea \be \sum_{i = 1}^N
\sum_{j = 1}^N \eta^{ij} \left ({\pa^2 \psi_{\alpha} \over \pa u^i
\pa u^q} {\pa^2 \psi_{\beta} \over \pa u^j \pa u^l} - {\pa^2
\psi_{\alpha} \over \pa u^i \pa u^l} {\pa^2 \psi_{\beta} \over \pa
u^j \pa u^q} \right ) = 0. \label{r3a} \ee

We consider a special {\it $k$-potential} ansatz for the functions
$\psi_{\alpha} (u)$, $1 \leq \alpha \leq L$. Consider an arbitrary
function $\Phi (u)$ and define the functions $\psi_{\alpha} (u)$, $1
\leq \alpha \leq L$, as follows: \be \psi_{(s - 1) N + i} = {\pa
\Phi \over \pa u^i}, \ \ \ \ 1\leq s \leq k, \ 1 \leq i \leq N,\ee
and $\psi_{\alpha} (u)$, $k N + 1 \leq \alpha \leq L,$ are arbitrary
functions that are linear in the coordinates (the corresponding
second fundamental forms vanish). In this case, the Gauss equations
(\ref{g3a}) can be rewritten in the form \bea && \sum_{r, s = 1}^k
\sum_{\alpha = (r - 1) N + 1}^{r N} \sum_{\beta = (s - 1) N + 1}^{s
N} \mu^{\alpha \beta} \left ({\pa^3 \Phi \over \pa u^{\alpha - (r -
1) N }\pa u^i \pa u^q} {\pa^3 \Phi \over \pa u^{\beta - (s - 1) N}
\pa u^j \pa u^l} - \right. \nn\\ && - \left.
 {\pa^3 \Phi \over \pa u^{\alpha - (r - 1) N}\pa u^i \pa u^l} {\pa^3
\Phi \over \pa u^{\beta - (s - 1) N } \pa u^j \pa u^q} \right ) = 0
\label{g3b} \eea or \be \sum_{r, s = 1}^k \sum_{m, n = 1}^N \mu^{(r
- 1) N + m, (s - 1) N + n} \left ({\pa^3 \Phi \over \pa u^m \pa u^i
\pa u^q} {\pa^3 \Phi \over \pa u^n \pa u^j \pa u^l} - {\pa^3 \Phi
\over \pa u^m \pa u^i \pa u^l} {\pa^3 \Phi \over \pa u^n \pa u^j \pa
u^q} \right ) = 0. \label{g3c} \ee

The Ricci equations (\ref{r3a}) in this case take the form \be
\sum_{i = 1}^N \sum_{j = 1}^N \eta^{ij} \left ({\pa^3 \Phi \over \pa
u^m \pa u^i \pa u^q} {\pa^3 \Phi \over \pa u^n \pa u^j \pa u^l} -
{\pa^3 \Phi \over \pa u^m \pa u^i \pa u^l} {\pa^3 \Phi \over \pa u^n
\pa u^j \pa u^q} \right ) = 0; \label{r3b} \ee i.e., they coincide
with the associativity equations of two-dimensional topological
quantum field theories (\ref{ass1a}) (the WDVV equations, see
\cite{1d}--\cite{4d}).

We consider now a special ansatz for the constant Gram matrices
$\mu_{\alpha \beta}$, $1 \leq \alpha, \beta \leq L$, in the normal
spaces of the submanifolds. Consider the case when \be \mu^{(r - 1)
N + m, (s - 1) N + n} = c^{rs} \eta^{mn}, \ \ \ \ 1 \leq r, s \leq
k, \ 1 \leq m, n \leq N, \label{gram} \ee where $c^{rs}$ is an
arbitrary nondegenerate symmetric constant matrix: $c^{rs} =
c^{sr}$, $\det(c^{rs}) \neq 0,$ $c^{rs} = {\rm \ const}$, and the
other elements of the matrix $\mu^{\alpha \beta}$ (for $\alpha \geq
k N + 1$ or $\beta \geq k N + 1$) are arbitrary constants such that
the matrix $\mu^{\alpha \beta}$ is symmetric and nondegenerate. For
such special constant Gram matrices $\mu_{\alpha \beta}$ (see
(\ref{gram})) the Gauss equations (\ref{g3c}) take the form \be
\sum_{r, s = 1}^k \sum_{m, n = 1}^N c^{rs} \eta^{mn} \left ({\pa^3
\Phi \over \pa u^m \pa u^i \pa u^q} {\pa^3 \Phi \over \pa u^n \pa
u^j \pa u^l} - {\pa^3 \Phi \over \pa u^m \pa u^i \pa u^l} {\pa^3
\Phi \over \pa u^n \pa u^j \pa u^q} \right ) = 0; \label{g3d} \ee
i.e., in this case, the Gauss equations (\ref{g3d}) are a linear
combination of the Ricci equations (\ref{r3b}) (the associativity
equations (\ref{ass1a})). Thus, in this case, all the fundamental
relations and equations of submanifold theory reduce to the
associativity equations (\ref{ass1a}). We will call such special
reductions of the fundamental Gauss--Codazzi--Ricci equations and
relations of submanifold theory {\it $k$-potential}.

\bt The associativity equations of two-dimensional topological
quantum field theories {\rm (\ref{ass1a})} are natural $k$-potential
reductions of the fundamental equations of submanifold theory.\et

By the Bonnet theorem, for any nondegenerate symmetric constant
matrix $c^{rs}$, $1 \leq r, s \leq k$, $c^{rs} = c^{sr}$,
$\det(c^{rs}) \neq 0,$ $c^{rs} = {\rm \ const}$, any nondegenerate
symmetric constant matrix $\eta^{ij}$, $1 \leq i, j \leq N$,
$\eta^{ij} = \eta^{ji},$ $\eta^{ij} = {\rm const},$ $\det
(\eta^{ij}) \neq 0$, and any nondegenerate symmetric constant matrix
$\mu^{\alpha \beta}$, $1 \leq \alpha, \beta \leq L$, $\mu^{\alpha
\beta}= {\rm const}$, $\mu^{\alpha \beta} = \mu^{\beta \alpha}$,
$\det (\mu^{\alpha \beta}) \neq 0$, such that relations (\ref{gram})
hold, any solution $\Phi (u)$ of the associativity equations
(\ref{ass1a}) that is determined up to quadratic terms gives a
unique (up to motion in the ambient pseudo-Euclidean space) totally
nonisotropic $N$-dimensional flat torsionless submanifold with the
first fundamental form $d s^2 = \eta_{ij} du^i du^j$, the second
fundamental forms $$\omega_{(s - 1) N + m, ij} (u) d u^i d u^j =
{\pa^3 \Phi \over \pa u^m \pa u^i \pa u^j}d u^i d u^j, \ \ \ \ 1
\leq s \leq k, \ 1 \leq i, j, m \leq N,$$
$$\omega_{p\, ij} (u) d u^i d u^j = 0, \ \ \ \ k N + 1 \leq p \leq L,
\ 1 \leq i, j \leq N,$$ and the constant Gram matrix $\mu_{\alpha
\beta}$, $1 \leq \alpha, \beta \leq L$, of the scalar products of
basic vectors in the normal spaces in an $(N + L)$-dimensional
pseudo-Euclidean space whose signature is the sum of the signatures
of the metrics $\eta_{ij}$, $1 \leq i, j \leq N$, and $\mu_{\alpha
\beta}$, $1 \leq \alpha, \beta \leq L$.

 We will call such submanifolds parametrized by the special
constant Gram mat\-ri\-ces $\mu_{\alpha \beta}$, $1 \leq \alpha,
\beta \leq L$ (see (\ref{gram})), and solutions of the associativity
equations (\ref{ass1a}) {\it $k$-potential}.

\bt The class of $k$-potential submanifolds in any Euclidean or
pseudo-Euclidean space is integrable. \et

\bt On each $k$-potential submanifold in a pseudo-Euclidean space,
there are $k$ natural identical structures of Frobenius algebras
{\rm (}$k$ identical Dubrovin--Frobeni\-us structures{\rm )} given
{\rm (}in flat coordinates{\rm )} for each $s$, $1 \leq s \leq k$,
by the first funda\-mental form $\eta_{ij}$ and the Weingarten
operators $(A_{(s - 1) N + m})^i_j (u) = - \eta^{il} \omega_{(s - 1)
N + m, lj} (u)$, $1 \leq i, j, l, m \leq N${\rm :} \bea && \langle
e_i, e_j \rangle = \eta_{ij}, \ \ \ \ e_i \circ e_j = c^l_{ij} (u)
e_l, \ \ \ \
e_i = {\pa \over \pa u^i}, \nn\\
&& c^l_{mj} (u^1, \ldots, u^N) = - (A_{(s - 1) N + m})^l_j (u) =
\eta^{li} \omega_{(s - 1) N + m, ij} (u^1, \ldots, u^N), \eea where
$\omega_{n\,ij} (u) d u^i d u^j,$ $1 \leq n \leq k N,$ are the
second fundamental forms of the submani\-fold. \et

\bt Each $N$-dimensional Frobenius manifold can be locally
represented as a $k$-potential $N$-dimensional submanifold in a $((k
+ 1) N + p)$-dimensional pseudo-Euclidean space {\rm (}for an
arbitrary positive integer $k$ and an arbitrary nonnegative integer
$p${\rm )}. \et

We note that the set of admissible signatures of the ambient
pseudo-Euclidean space can be easily determined by the signature of
the metric $\eta_{ij}$, $1 \leq i, j \leq N$, of the Frobenius
manifold and by the given integers $k$ and $p$ (this set is never
empty). Let $2s - N$ be the signature of the metric $\eta_{ij}$ of
the Frobenius manifold, where $s$, $0 \leq s \leq N$, is the
positive index of inertia of the metric. Then the set of admissible
signatures of the ambient pseudo-Euclidean space is determined by
the formula $(2 s - N) (2 r - k + 1) + 2t - p$, $0 \leq r \leq k$,
$0 \leq t \leq p$. In particular, in the simplest case when $p = 0$
and $k = 1$, only the signatures $2 (2s - N)$ and $0$ are
admissible.

\bt For an arbitrary Frobenius manifold that is given locally by a
solution $\Phi (u^1, \ldots, u^N)$ of the associativity equations
{\rm (\ref{ass1a})}, the corresponding $k$-potential subma\-nifolds
in $((k + 1) N + p)$-dimensional pseudo-Euclidean spaces that
realize this Frobenius manifold are determined by any $((k + 1) N +
p)$-component vector function $r (u^1, \ldots, u^N)$ satisfying the
following consistent linear system of second-order partial
differential equations{\rm :} \bea && {\partial^2 r \over
\partial u^i \partial u^j} = \sum_{r, s = 1}^k \sum_{m, l = 1}^N c^{rs}
\eta^{ml} {\partial^3 \Phi \over
\partial u^i \partial u^j \partial u^m} n_{(s - 1) N + l} (u) +\nn\\
&& + \sum_{\beta = k N + 1}^{k N + p} \sum_{r = 1}^k \sum_{m = 1}^N
\mu^{(r - 1 ) N + m, \beta} {\partial^3 \Phi \over
\partial u^i \partial u^j \partial u^m} n_{\beta} (u),\ \ \ \ 1 \leq i, j \leq N,
 \label{rf1} \eea \be {\partial n_{(r - 1) N + m} \over \partial u^i}
 = - \sum_{l, j = 1}^N \eta^{jl} {\partial^3 \Phi \over \partial u^l \partial u^m
\partial u^i} {\partial r \over \partial u^j},
 \ \ \ \ 1 \leq r \leq k, \ 1 \leq i, m \leq N, \label{rf2} \ee
 \be {\partial n_{\alpha} \over \partial u^i}
 = 0,
 \ \ \ \ k N + 1 \leq \alpha \leq k N + p, \ 1 \leq i \leq N, \label{rf2a} \ee
where $n_{\alpha} (u^1, \ldots, u^N)$, $1 \leq \alpha \leq k N + p,$
are some $((k + 1) N + p)$-component vector functions. The
consistency conditions for the linear system {\rm
(\ref{rf1})--(\ref{rf2a})} are equivalent to the associativity
equations {\rm (\ref{ass1a})} for the function $\Phi (u)$. \et

{\bf Acknowledgments.} The work was supported by the
Max-Planck-Institut f\"ur Mathematik (Bonn, Germany), by the Russian
Foundation for Basic Research (project no.~08-01-00054) and by a
grant of the President of the Russian Federation (project no.
NSh-1824.2008.1).

\begin{flushleft}
{\bf O. I. Mokhov}\\
Centre for Nonlinear Studies,\\
L.D.Landau Institute for Theoretical Physics,\\
Russian Academy of Sciences,\\
Kosygina str., 2, Moscow, Russia;\\
Department of Geometry and Topology,\\
Faculty of Mechanics and Mathematics,\\
M.V.Lomonosov Moscow State University,\\
Moscow, Russia\\
{\it E-mail\,}: mokhov@mi.ras.ru; mokhov@landau.ac.ru; mokhov@bk.ru\\
\end{flushleft}

\end{document}